\documentclass[12pt]{article}
\newcommand{\Hstar}{H^{*}}
\newcommand{\ract}{\mathbin{{\leftharpoonup}}}
\newcommand{\lact}{\mathbin{{\rightharpoonup}}}
\newcommand{\Ker}{\mathop{\rm Ker}}
\renewcommand{\Im}{\mathop{\rm Im}}
\newcommand{\End}[1]{\mathop{{\rm End}_{#1}}}
\newcommand{\tensor}[1]{\mathbin{\otimes_{#1}}}
\newcommand{\ebar}{{\bar e}}
\newcommand{\xbar}{{\bar x}}
\newcommand{\xibar}{{\bar \xi}}

\newcommand{\eu}[1]{\mbox{$ {\bf R}^{#1}$}}

\newcommand{\jperp}{\mbox{$ J^{\perp} $}}

\newcommand{\LT}{{\cal LT}}

\newcommand{\assoc}[3]{\mbox{${#3} {<} {#1}_1, \ldots, {#1}_{#2} {>}$}}
\newcommand{\prtl}[1]{\mathop{\frac{\partial}{\partial x_{#1}}}}
 
\newtheorem{theorem}{Theorem}[section]
\newtheorem{lemma}{Lemma}[section]
\newtheorem{corollary}{Corollary}[section]
\newcommand{\proof}{\medskip\noindent{\sc Proof.\,}}
\newcommand{\closeproof}{\hfill\hbox{\vrule width1.0ex height1.5ex}\medskip}
 
\title{The realization of input-output maps using
bialgebras}
\author{Robert Grossman\thanks{This research was
supported in part by the grants
NASA NAG2-513 and NSF DMS-8904740.} \,
and Richard G. Larson\thanks{Supported in
part by NSF Grant DMS 870--1085.}\\
University of Illinois at Chicago}
\date{February, 1991}
\begin{document}
\maketitle
 
\medbreak
{\noindent\bf
This is a draft of a paper which later appeared 
in Forum Mathematicum, Volume 4, pp. 109-121, 1992.}
\medbreak

\begin{abstract}
We use the theory of bialgebras to provide the
algebraic background for state space
realization theorems for input-output maps
of control systems.  This allows us to consider
from a common viewpoint classical results about
formal state space realizations of nonlinear systems
and more recent results involving analysis related to
families of trees.
If $H$ is a bialgebra, we say that
$p\in\Hstar$ is {\em differentially produced by the algebra $R$ with
the augmentation $\epsilon$\/} if
there is right $H$-module algebra structure on $R$
and there exists $f\in R$ satisfying
$p(h) = \epsilon(f\cdot h)$.
We characterize those $p\in\Hstar$ which are differentially produced.
\end{abstract}

\section{Introduction}\label{section:introduction}
In this paper, we use the theory of bialgebras to provide the
algebraic background for state space
realization theorems for input-output maps
of control systems.  This approach allows us to consider
from a common viewpoint the formal aspects
of classical results of
Fliess~\cite{Fliess:realization1},~\cite{Fliess:realization2}
and more recent results involving analysis related to
families of trees \cite{GL:solve},~\cite{GL:families}.
The following definition is fundamental to this approach.
If $H$ is a bialgebra, we say that
$p\in\Hstar$ is {\em differentially produced by the algebra $R$ with
the augmentation $\epsilon$\/} if
\begin{enumerate}
\item there is right $H$-module algebra structure on $R$
(see the definition of $H$-module algebra below);
\item there exists $f\in R$ satisfying
$p(h) = \epsilon(f\cdot h)$.
\end{enumerate}
We will characterize those $p\in\Hstar$
which are differentially produced.
 
Differentially produced elements of algebras arise naturally
when studying control systems with inputs and outputs.
For example, let $R$ denote the field of rational functions in the
variables
$x_1$, \ldots, $x_N$ with coefficients from the
field $k$, and let  $E_1$, \ldots, $E_M$ denote $M$ derivations of
$R$. The control system
\begin{eqnarray}\label{ControlSystem}
\dot{x}(t) &=& \sum_{i=1}^M u_{i}(t) E_{i}(x(t)), \nonumber\\
x(0) &=& x^0\in\eu{N}\end{eqnarray}
together with an observation function $f\in R$
\begin{equation} f:\eu{N}\longrightarrow\eu{}
\label{ObservationFunction}\end{equation}
naturally specifies an input-output map, which is
defined by sending the input functions
$$ t\rightarrow u_1(t), \quad\ldots, \quad t\rightarrow u_M(t)$$
to the output function
$$ t\rightarrow f(x(t)).$$
The properties of the input-output map are captured
by the formal series
$$\sum_{\mbox{\scriptsize words } \mu} c_{\mu} {\mu},$$
where
$$ c_{\mu}=E_{\mu_k} \cdots E_{\mu_1} f(x(0)),\qquad
\mbox{$\mu=\mu_1\cdots\mu_k$ a word.}$$
This series is often called the generating series,
while the data consisting of a control system with inputs,
together with an observation, are called a state space
realization of the input-output map.
Isidori \cite{Isidori} contains a detailed description
of these topics, as well as extensive references.
 
The results in this paper do not address the analytical aspects
of the realization of input-output maps or questions of convergence
of the series involved, but are purely algebraic and combinatorial.
 
Let $H$ denote the free associative algebra in
the symbols $E_1,\ldots, E_M$ over the
field $k$ and let $\Hstar$ denote its topological dual.
$\Hstar$ is isomorphic to a formal power series algebra
in infinitely many variables.
The point of view of this paper is to consider the formal series
$p$ as an element of the algebra $\Hstar$.
If $p\in\Hstar$ is the formal series associated with
an input-output map, then it is differentially produced.
Conversely, we can ask which formal series $p\in\Hstar$
have the property that there is a control system and
an observation function which realizes it as above; that
is, which $p$ are differentially produced?  We will see
that both these questions are simply answered if we exploit
the bialgebra structure of $H$.
 
Important work in this area has been done
by Fliess \cite{Fliess:realization1},~\cite{Fliess:realization2},
Hermann and Krener~\cite{HermannKrener}, and Sussman~\cite{Sussman}.
Fliess was the first to focus
on the algebraic and combinatorial
aspects of the problem, making important use of
shuffle algebras in his study of realization theory.  His work
was simplified by Reutenauer~\cite{Reutenauer}.
In this paper, we generalize and simplify the
work of Fliess  and Reutenauer, extending the context
to general bialgebras.  This also allows us to treat
combinatorial examples of differential representations
which have arisen in the symbolic computation of
solutions of differential equations~\cite{GL:effecient1}
and~\cite{GL:solve}.
An alternate combinatorial treatment of generating series is due to
Viennot~\cite{viennot}.
Hoang Ngoc Minh's work \cite{Minh} on the evaluation transform is also
related to the work described here.
See the survey article by Jakubczyk~\cite{Jakubczyk} for other
approaches to the realization theory of nonlinear systems.
 
To state the theorem we prove in Section~\ref{section:proof},
we need some definitions.
Let $k$ denote a field of characteristic $0$.
If $V$ is a vector space over $k$, denote by $V^{*}$ the set of all
linear maps $V\rightarrow k$.
The vector space $V^{*}$ with the finite topology
is a complete topological vector space (see~\cite{NJ:StrRings} for
details).
 
By an {\em algebra\/} over $k$ we mean an associative algebra with
identity.
The algebra structure of $A$ can be specified by the maps
$A\tensor{k}A\rightarrow A$ which maps $a\otimes b\in A\tensor{k}A$ to
$ab\in A$, and $k\rightarrow A$ which maps $1\in k$ to $1\in A$.
The facts that multiplication is associative and that $1\in A$ is a
two-sided unit for multiplication can be expressed in terms of the
commutativity of certain diagrams.
An {\em augmentation\/} for the algebra $A$ is an algebra homomorphism
$A\rightarrow k$.
 
In a dual manner we define a {\em coalgebra\/}.
A coalgebra is a vector space $C$ over $k$, equipped with maps
$\Delta : C\rightarrow C\tensor{k}C$ and $\epsilon : C\rightarrow k$
which give a coassociative comultiplication and a counit
(that is, the diagrams which are dual to those in the definition of an
algebra commute).
If $c\in C$ we will sometimes write the element $\Delta(c)\in
C\tensor{k}C$ as $\sum_{(c)}c_{(1)}\otimes c_{(2)}$, using notation
introduced by Sweedler in~\cite{Moss}.
If $C$ is a coalgebra, then $C^{*}$ is an algebra.
Note that $\epsilon : C\rightarrow k$ is the multiplicative identity for
the algebra $C^{*}$.
 
We define a {\em bialgebra\/} to be a vector space $H$ equipped with an
algebra and a coalgebra structure, so that the maps which define the
coalgebra structure are algebra homomorphisms, or equivalently, the maps
which define the algebra structure are coalgebra homomorphisms.
In particular, the coalgebra counit $\epsilon$ is an augmentation for
the algebra $H$.
If $H$ is a bialgebra, its {\em primitive\/} elements are defined by
\[
P(H) = \{\,h\in H \mid \Delta(h)=1\otimes h + h\otimes1\,\}.
\]
It can be shown that $P(H)$ is a Lie algebra with respect to the
operation $[x,y]=xy-yx$.
It has been shown~\cite{MM} that if the characteristic of~$k$ is~0,
and $H$ is generated as
an algebra by $P(H)$ (in which case we say that $H$ is {\em primitively
generated\/}), then $H\cong U(P(H))$, where $U(L)$ denotes the universal
enveloping algebra of the Lie algebra $L$.
The Poincar\'e-Birkhoff-Witt Theorem (see~\cite{NJ:LieAlgs}) states that
if
$e_1$, $e_2$,~\ldots\ is an ordered basis for $L$, then
\[
\{\,e_{i_1}^{\alpha_{i_1}} \cdots e_{i_k}^{\alpha_{i_k}} \mid
i_1 < \cdots < i_k\mbox{ and } 0<\alpha_{i_r}\,\}
\]
is a basis for $U(L)$.
It follows that $U(L)^{*}$ is a formal power series algebra.
More specifically, if we denote the basis element
$e_{i_1}^{\alpha_{i_1}} \cdots e_{i_k}^{\alpha_{i_k}}$ of $U(L)$ by
$e^{\alpha}$,
and let $\{x_{\alpha}\}$ be the dual basis (in the sense of complete
topological vector spaces), then
$U(L)^{*}\cong k[[x_1,x_2,\ldots\,]]$, where $x_i$ is the dual basis
element corresponding to $e_i$.
Under this isomorphism
\[
x_{\alpha}=\frac{x^{\alpha}}{\alpha!},
\]
where $x^{\alpha}=x_{i_1}^{\alpha_{i_1}} \cdots x_{i_k}^{\alpha_{i_k}}$
and $\alpha!=\alpha_{i_1}!\cdots\alpha_{i_k}!$.
(Note that we can think of $\alpha$ as an infinite sequence of
non-negative integers, all but finitely many of which are 0.
Recall that $x_j^0=1$ and $0!=1$.)
 
 
Let $H$ be a primitively generated bialgebra.
We define a right and left $H$-module structure on $\Hstar$ as follows:
if $p\in\Hstar$ and $h\in H$, let $p\ract h\in\Hstar$ be defined by
\[
 (p \ract h)(k) = p(hk), \qquad k\in H,
\]
and let $h\lact p\in\Hstar$ be defined by
\[
 (h \lact p)(k) = p(kh), \qquad k\in H.
\]
We say that an algebra $A$ is a {\em left $H$-module algebra\/} if $A$
is a left $H$-module, and
\[
h\cdot(ab)=\sum_{(h)}(h_{(1)}\cdot a)(h_{(2)}\cdot b).
\]
A {\em right $H$-module algebra\/} is defined similarly.
If $A$ is a left or right $H$-module algebra, we say that {\em $H$
measures $A$ to itself\/}.
In particular, $H$ measures $\Hstar$ to itself using the actions $\lact$
and $\ract$ defined above.
If the bialgebra $H$ measures the algebra $A$ to itself, then the
elements of $P(H)$ act as derivations of $A$.
 
We say that $p\in\Hstar$ has {\em finite
Lie rank\/} if $\dim P(H)\lact p$ is finite.
Recall that
$p\in\Hstar$ is {\em differentially produced by the algebra $R$ with
the augmentation $\epsilon$\/} if
\begin{enumerate}
\item there is right $H$-module algebra structure on $R$;
\item there exists $f\in R$ satisfying
$p(h) = \epsilon(f\cdot h)$.
\end{enumerate}
Concrete examples of differentially produced functionals on a
primitively generated bialgebra (that
is, of differentially produced formal power series) are given in
Section~\ref{section:examples}.
 
Our main theorem is the following.
\begin{theorem}\label{maintheorem}
Let $H$ be a primitively generated bialgebra
over a field of characteristic $0$.
Let $p\in\Hstar$.
Then the following are equivalent:
\begin{enumerate}
\item\label{ThmItem1} $p$ has finite Lie rank;
\item\label{ThmItem2} $p$ is differentially produced by some augmented
$k$-algebra
$R$ for which $\dim\,(\Ker\epsilon)/(\Ker\epsilon)^2$ is finite;
\item\label{ThmItem3} $p$ is differentially produced by a subalgebra of
$\Hstar$ which
is isomorphic to $k[[x_1,\ldots,x_N]]$, the algebra of formal power
series in $N$ variables.
\end{enumerate}
\end{theorem}
 
We prove this theorem in Section~\ref{section:proof}.
We give examples of its application in Section~\ref{section:examples}.
 
\section{Proof of Main Theorem}\label{section:proof}
We first prove that part~(\ref{ThmItem1}) of
Theorem~\ref{maintheorem} implies part~(\ref{ThmItem3}).
Given a fixed $p\in\Hstar$, we define three basic objects:
\begin{eqnarray*}
L &=& \{\, h \in P(H) \mid h\lact p = 0\,\} \\
J &=& HL \\
\jperp &=& \{\,q\in\Hstar \mid q(j)=0 \mbox{ for all } j \in J \}.
\end{eqnarray*}
Since $L\subseteq P(H)$, it follows that $J$ is a coideal,
that is, that $\Delta(J)\subseteq J\otimes H + H\otimes J$.
Therefore $\jperp\cong(H/J)^{*}$ is a subalgebra of $\Hstar$.
We will show that $\jperp$ is
isomorphic to a formal power series algebra, and
will construct derivations of this ring which will be used
to realize the input-output map defined by $p$.
 
\begin{lemma}\label{lemma-a}
If
$\dim P(H)\lact p = N$,
then $\jperp$ is a subalgebra of $\Hstar$ satisfying
\[
\jperp \cong k[[x_1,\ldots, x_N]].
\]
\end{lemma}
 
\proof
The sub Lie algebra $L$ has finite codimension $N$.
Choose a basis $\{ e_1, e_2, \ldots \}$ of $P(H)$ such
that $\{ e_{N+1}, e_{N+2}, \ldots \}$ is a basis of $L$.
Note that if $\ebar_i$ is
the image of $e_i$
under the quotient map
$ P(H) \rightarrow P(H)/L$, then
$\{ \ebar_1$, \ldots, $\ebar_N \}$
is a basis for $P(H)/L$.
 
By the Poincar\'e-Birkhoff-Witt Theorem, $H$ has
a basis of the form
\[
\{\,e_{i_1}^{\alpha_{i_1}} \cdots e_{i_k}^{\alpha_{i_k}} \mid
i_1 < \cdots < i_k\mbox{ and } 0<\alpha_{i_r}\,\}.
\]
Since $L$ is a sub Lie algebra of $P(H)$, and the basis $\{e_i\}$ of
$P(H)$ has been chosen so that $e_i\in L$ for $i>N$, it follows that the
operation of putting monomials in standard form which is used in the
proof of the Poincar\'e-Birkhoff-Witt Theorem will map elements of
$J=HL$ to linear combinations of monomials
of the form
\[
e_{i_1}^{\alpha_{i_1}} \cdots e_{i_k}^{\alpha_{i_k}}
\]
with at least one $i_r>N$.
Therefore $J$ has a basis of such monomials.
It follows that
\[
\{\,\ebar_1^{\alpha_1}\cdots\ebar_N^{\alpha_N} \mid
\alpha_1,\ldots,\alpha_N\ge0\,\}
\]
is a basis for $H/J$.
It now follows that the elements of the form
\[
x_{\alpha} =
\frac{x^{\alpha}}{\alpha !} =
\frac{ x_{i_1}^{\alpha_{i_1}} \cdots x_{i_k}^{\alpha_{i_k}}}
{\alpha_{i_1}! \cdots \alpha_{i_k}! }
\]
with all $1\le i_r\le N$ are in $\jperp\subseteq\Hstar$.
Indeed, $\jperp$ consists
precisely of the completion in the finite topology of the span of such
elements. In other words,
\[
\jperp \cong k[[x_1,\ldots, x_N]],
\]
completing the proof.
\closeproof
 
We will use the following notation and facts from the proof of
Lemma~\ref{lemma-a}:
Suppose that $\{e_1$, \ldots, $e_N$, $\ldots\}$ is a basis for $P(H)$
such that
$\{e_{N+1}$, \ldots $\}$ is a basis for $L$.
Let $\{e^{\alpha}\}$ be the corresponding Poincar\'e-Birkhoff-Witt
basis.
Denote $\jperp$ by $R$.
Then $R\cong k[[x_1$, \ldots, $x_N]]$, and
$x_1^{\alpha_1}\cdots x_N^{\alpha_N}/{\alpha_1}!\cdots{\alpha_N}!$
equals the element of the dual (topological) basis of $\Hstar$
to the Poincar\'e-Birkhoff-Witt basis $\{e^{\alpha}\}$ of $H$,
corresponding
to the basis element $e_1^{\alpha_1}\cdots e_N^{\alpha_N}$.
 
We now collect some properties of the ring of formal power
series $R$ which will be necessary for the proof of the theorem.
 
\begin{lemma}\label{lemma-b}
Assume $p\in\Hstar$ has finite Lie rank, and
let $R\subseteq \Hstar$, $e_{\alpha}\in H$, and $x^{\alpha}\in R$ be as
above. Define
\[
f = \sum_{\scriptstyle\alpha=(\alpha_1,\ldots,\alpha_N)\atop
\scriptstyle\alpha\geq 0} c_{\alpha}x^{\alpha}\in R,
\]
where $c_{\alpha} = \frac{\textstyle p(e^{\alpha})}{\textstyle\alpha!}$.
Then
\begin{enumerate}
\item\label{AssMeas} $H$ measures $R$ to itself via $\ract$;
\item\label{AssReal} $p(h)=\epsilon(f\ract h)$ for all $h\in H$.
\end{enumerate}
\end{lemma}
 
\proof We begin with the proof of part~(\ref{AssMeas}).
Since $H$ measures $\Hstar$ to itself and $R\subseteq\Hstar$,
we need show only that $R\ract H\subseteq R$.
Take $r\in R$, $h\in H$ and $j\in J$.
We have $ (r\ract h)(j) = r(hj)$.
Since $J$ is a left ideal, $hj\in J$, so $r(hj)=0$, so
$r\ract h \in\jperp = R$.
This proves part~(\ref{AssMeas}).
 
We now prove part~(\ref{AssReal}).
Let $e^{\alpha}=e_{i_1}^{\alpha_{i_1}} \cdots e_{i_k}^{\alpha_{i_k}}$
be a Poincar\'e-Birkhoff-Witt basis element of $H$.
Since $e^{\alpha}\in J$
unless $\{i_1,\ldots,i_k\}\subseteq\{1,\ldots,N\}$,
$p(e^{\alpha})=0$
unless $\{i_1,\ldots,i_k\}\subseteq\{1,\ldots,N\}$.
Also $\epsilon(f\ract e^{\alpha})=f\ract e^{\alpha}(1)=
f(e^{\alpha}1)=f(e^{\alpha})=0$
unless $\{i_1,\ldots,i_k\}\subseteq\{1,\ldots,N\}$.
Now suppose $\{i_1$, \ldots, $i_k\}\subseteq\{1,\ldots,N\}$.
We have in this case that
$p(e^{\alpha})=\alpha!c_{\alpha}=f(e^{\alpha})=
f\ract e^{\alpha}(1)=\epsilon(f\ract e^{\alpha})$.
Since $\{e^{\alpha}\}$ is a basis for $H$, this completes the proof
of part~(\ref{AssReal}) of the lemma.
\closeproof
 
\begin{corollary}
Under the assumptions of Lemma~\ref{lemma-b}, $f=p$.
\end{corollary}
 
Lemmas~\ref{lemma-a} and~\ref{lemma-b} yield that
part~(\ref{ThmItem1}) implies part~(\ref{ThmItem3}) in
Theorem~\ref{maintheorem}.
It is immediate that
part~(\ref{ThmItem3}) implies part~(\ref{ThmItem2}).
 
We now complete the proof of Theorem~\ref{maintheorem} by proving that
part~(\ref{ThmItem2}) implies part~(\ref{ThmItem1}).
 
Let $x_1,\ldots,x_N$ be chosen so that
$\{\xbar_1,\ldots,\xbar_N\}$ is a basis for
$(\Ker\epsilon)/(\Ker\epsilon)^2$.
If $f\in R$ and $h\in H$, then
\[
f\cdot h=q_0(h)1+\sum_{i=1}^{N}q_i(h)x_i+g(h),
\]
where $q_i\in\Hstar$ and $g(h)\in (\Ker\epsilon)^2$.
Let $l\in P(H)$.
Since $H$ measures $R$ to itself and $\Delta(l)=1\otimes l+l\otimes1$,
the map $f\mapsto f\cdot l$ is a derivation of $R$.
Now let $f\in R$ be the element such that
\[
p(h)=\epsilon(f\cdot h).
\]
Then
\begin{eqnarray*}
f\cdot hl & = & (f\cdot h)\cdot l \\
          & = & q_0(h)1\cdot l+\sum_{i=1}^N q_i(h)x_i\cdot l + g(h)\cdot l.
\end{eqnarray*}
Since the map $f\mapsto f\cdot l$ is a derivation,
$1\cdot l=0$;
since $g(h)\in (\Ker\epsilon)^2$, $g(h)\cdot l\in\Ker\epsilon$.
It follows that
\begin{eqnarray*}
l\lact p(h) & = & p(hl)\\
            & = & \epsilon(f\cdot hl)\\
            & = & \sum_{i=1}^Nq_i(h)\epsilon(x_i\cdot l).
\end{eqnarray*}
Therefore $P(H)\lact p\subseteq \sum_{i=1}^Nkq_i$, so $p$ has finite Lie
rank.
This completes the proof of Theorem~\ref{maintheorem}
 
\section{Examples}\label{section:examples}
In this section, we discuss two examples
of applications of Theorem~\ref{maintheorem}.
The first example is obtained by letting the bialgebra
$H$ be the
free associative algebra over the field $k$ in the symbols
$E_1,\ldots, E_M$.  This example motivated the theorem and is
the basic setting in the control theory literature
(see~\cite{Fliess:realization1}, \cite{Fliess:realization2},
and \cite{Reutenauer}).  The second example is
obtained by letting the  bialgebra $H$ have as basis some family
of trees with the appropriate
multiplication and comultiplication.  This example arises
when studying algorithms for the symbolic computation
of higher order derivations generated by derivations
$E_1,\ldots, E_M$; see~\cite{GL:effecient1} and~\cite{GL:effecient2}.
There is a natural homomorphism between
these two Hopf algebras which is described in \cite{GL:effecient1}.
 
\paragraph{\bf Example 1.}
We begin by giving a description of the setting
for this example.
Let $R$ denote the field of rational functions in the variables
$x_1,\ldots, x_N$ with coefficients from the
field $k$, and let  $E_1,\ldots, E_M$ denote $M$ derivations of
$R$. The algebras in this example are associated with
a pair consisting of the dynamical system~(\ref{ControlSystem})
and the observation function~(\ref{ObservationFunction})
introduced in Section~\ref{section:introduction}.
We assume that the controls
$$t\rightarrow u_1(t),\ldots, \quad t\rightarrow u_M(t)$$
in (\ref{ControlSystem}) are
continuous and bounded almost everywhere.
 
Integrating the initial value problem (\ref{ControlSystem}) gives
\[
f(x(t)) = f(x(0)) + \sum_{\mu_1=1}^M \int_{0}^{t}
u_{\mu_1}(\tau) E_{\mu_1}(f(x(\tau)))\, d\tau.
\]
Integrating again gives
\begin{eqnarray*}
f(x(t)) &=& f(x(0)) + \sum_{\mu_1=1}^M
E_{\mu_1}f(x(0))
\int_{0}^{t} u_{\mu_{1}}(\tau) \, d\tau_1  \\
&+& \sum_{\mu_1,\mu_2 =1}^M \int_{0}^{t} \int_{0}^{\tau_1}
u_{\mu_{1}}(\tau_1) u_{\mu_{2}} (\tau_2)
E_{\mu_1} E_{\mu_2} (f(x(\tau_2)))\, d\tau_2 d\tau_1.
\end{eqnarray*}
Continuing this process yields
\begin{eqnarray*}
f(x(t)) &=& f(x(0)) + \sum_{\mu_1=1}^M E_{\mu_1}(f(x(0)))
\int_{0}^{t} u_{\mu_{1}}(\tau) \, d\tau_1  \\*
&& {} + \cdots + \sum_{\mu_1,\ldots, \mu_{k} =1}^M
E_{\mu_{k}} \cdots E_{\mu_1} f(x(0))\cdot{} \\*
&& \qquad \int_{0}^{t} \cdots \int_{0}^{\tau_{k-1}}
u_{\mu_{1}}(\tau_1) \cdots u_{\mu_{k}}(\tau_{k}) \,
d\tau_{k} \cdots d\tau_1 + {\cal R},
\end{eqnarray*}
where the remainder $\cal R$ is of the form
\[
\sum_{\mu_1,\ldots,\mu_{k+1} =1}^M
\int_{0}^{t} \cdots \int_{0}^{\tau_{k}}
u_{\mu_{1}}(\tau_1) \cdots u_{\mu_{k+1}}(\tau_{k+1}) $$
$$ \qquad\qquad {} \cdot E_{\mu_{k+1}} \cdots E_{\mu_1}
f(x(\tau_{k+1}))\, d\tau_{k+1} \cdots d\tau_1.
\]
 
Let $\mu=\mu_1\cdots\mu_k$
denote a word of length $k$
built from the alphabet $\{1$, \ldots, $M\}$.
The above process defines a formal series
\begin{equation}\label{FormalSeries}
\bar{p}=\sum_{\mbox{\scriptsize words } \mu}
c_{\mu} {\bar{\xi}}_{\mu},
\end{equation}
where
\begin{equation} c_{\mu} = E_{\mu_k} \cdots E_{\mu_1}
f(x(0))\in k,\label{CEquations}\end{equation}
and
\begin{equation}\label{XiEquations}
{\bar{\xi}}_{\mu}(t) = \left\{
\begin{array}{ll}
\displaystyle
\int_{0}^t u_{\mu_1}(\tau)\, d\tau & \qquad\mbox{if $\mu=\mu_1$} \\
\displaystyle
\int_{0}^t u_{\mu_1}(\tau) {\bar{\xi}}_{\mu_2\cdots\mu_k}(\tau) d\tau
& \qquad\mbox{if $\mu=\mu_1 \cdots \mu_k$.}
\end{array}\right.
\end{equation}
Chen~\cite{Chen:integration} proved that functions
of the form (\ref{XiEquations}) form a shuffle algebra, that is, that
\[
\xibar_{\lambda}(t) \cdot \xibar_{\mu}(t)
= \sum_{\nu} \xibar_{\nu}(t),
\]
where the sum is over words $\nu$ that are in the shuffle
of the words $\lambda$ and $\mu$.
The shuffle of two words
\[
\alpha=\alpha_1\cdots\alpha_i, \qquad
\beta=\beta_1\cdots\beta_j
\]
is defined as follows.
Let
$K=\{1,2,\ldots, i+j\}$, and let
\begin{eqnarray*} \lambda &:& \{1,\ldots, i\} \longrightarrow K\\
\mu &:& \{1,\ldots, j\} \longrightarrow K
\end{eqnarray*}
denote two order-preserving maps such that the images
$\lambda(\alpha)$ and $\mu(\beta)$ are disjoint and
complementary.  These data define a word
$\gamma=\gamma_1\cdots\gamma_{k}$ via
\[
\gamma_l = \left\{
\begin{array}{rl}
\alpha_{\lambda^{-1}(l)} & \qquad\mbox{if $l\in\Im\lambda$;} \\
\beta_{\mu^{-1}(l)}      & \qquad\mbox{if $l\in\Im\mu$.}
\end{array}\right.
\]
The {\em shuffle} of $\alpha$ and $\beta$ is defined to be
the set of all such $\gamma$ obtained in this fashion.
We can now define a shuffle algebra structure.
Suppose a vector space has a spanning set $\{\xi_{\alpha}\}$ indexed by
all words $\alpha=\alpha_1\cdots\alpha_i$ over some alphabet.
The {\em shuffle product} of two elements
$\xi_{\alpha}$ and $\xi_{\beta}$ is defined by
$$ \xi_{\alpha} \cdot \xi_{\beta} =
\sum_{\gamma\in\Gamma} \xi_{\gamma},$$
where $\Gamma$ is the shuffle of the words
$\alpha$ and $\beta$.
The algebra $H$ is called a shuffle algebra if the multiplication in $H$
(with respect to some spanning set) is given by the shuffle
product.
More details on the shuffle algebra can be found in~\cite{Moss}.
 
Let $H=\assoc{E}{M}{k}$ denote the free
associative algebra in the symbols
$E_1,\ldots, E_M$ over the field $k$.
Recall that $H$ is a bialgebra.
The coproduct and counit are defined by letting
\begin{eqnarray*}
\Delta(E_i)   & = & 1\otimes E_i + E_i\otimes 1, \\
\epsilon(E_i) & = & 0,
\end{eqnarray*}
for $i=1$,\ldots, $M$, and then
extending to all of $\assoc{E}{M}{k}$
by requiring that $\Delta$ and $\epsilon$ be algebra homomorphisms.
The  bialgebra $H$ is cocommutative, but not commutative.
The algebra of formal series in the ${\bar{\xi}}_{\mu}(t)$
is a quotient of the algebra $\Hstar$.
 
The papers of
Fliess~\cite{Fliess:realization1}, \cite{Fliess:realization2},
Reutenauer \cite{Reutenauer}, and Crouch and Lamnabhi-Lagarrigue
\cite{CL:algebraic} all view the formal series $\bar{p}$ above
as an element of the shuffle algebra of formal power series
in the noncommuting variables $E_1,\ldots, E_M$.
It is easy to relate that point of view to the point of view taken here.
The bialgebra $H$ has basis consisting of all
words $E_{\mu}$ in the generators $E_1,\ldots, E_M$,
including the empty word $1$.  Let
$\xi_{\nu}$ denote the elements in the
dual $\Hstar$ of $H$ which are
formally dual to the $E_{\mu}$, that is,
$$ \xi_{\nu}(E_{\mu}) = \delta_{\mu, \nu},$$
where $\delta_{\mu, \nu}$ is the Kronecker delta.
Then the $\xi_{\nu}$ can be viewed as a topological basis for the
formal non-commutative power series ring over $k$ generated by
$E_1$, \ldots, $E_M$.
The algebra $\Hstar$ is a commutative algebra with
respect to the shuffle product on the $\xi_{\mu}$.
Fixing a control system~(\ref{ControlSystem}) and
an observation function~(\ref{ObservationFunction})
determines an element of $\Hstar$
\[
p= \sum_{\mbox{\scriptsize words } \mu} c_{\mu} \xi_{\mu},
\]
where the $c_{\mu}$ are given by Equation~(\ref{CEquations}).
The element
$\bar{p}$ given by Equation~(\ref{FormalSeries})
can be viewed as an element of a quotient
algebra of $\Hstar$.
Theorem~\ref{maintheorem} applied to this example gives
the classical theorem
of Fliess~\cite{Fliess:realization1},~\cite{Fliess:realization2}.
 
\paragraph{\bf Example 2.}
We follow \cite{GL:effecient1},~\cite{GL:effecient2} for this
example.
The algebra of trees we discuss here is important in the development of
efficient algorithms for the solution of differential equations.
By a tree we mean a finite rooted tree~\cite{Tarjan}.
If $\{E_1$, \ldots, $E_M\}$ is a set of symbols, we will say a tree
is {\em labeled with} $\{E_1$, \ldots, $E_M\}$ if every node of the
tree other than the root has an element of $\{E_1$, \ldots, $E_M\}$
assigned to it.
We denote the set of all trees
labeled with $\{E_1$, \ldots, $E_M\}$ by $\LT(E_1$, \ldots, $E_M)$.
Let
$k\{\LT(E_1$, \ldots, $E_M)\}$
denote the vector space over $k$ with basis $\LT(E_1$, \ldots, $E_M)$.
We show that this vector space is a graded connected Hopf
algebra.
 
We define the multiplication in $k\{\LT(E_1$, \ldots, $E_M)\}$ as
follows.
Since the set of labeled trees form a basis for $k\{\LT(E_1$, \ldots,
$E_M)\}$, it is
sufficient to describe the product of two labeled trees.
Suppose $t_1$ and $t_2$ are two labeled trees.
Let $s_1$, \ldots, $s_r$ be the children of the root of $t_1$.
If $t_2$ has $n+1$ nodes (counting the root), there are
$(n+1)^r$ ways to attach the $r$ subtrees of $t_1$
which have $s_1$, \ldots, $s_r$ as roots to the labeled tree $t_2$
by making each $s_i$ the child of some node of $t_2$,
keeping the original labels.
The product $t_1t_2$ is defined to be the sum of these $(n+1)^r$
labeled trees.
It can be shown that this product is associative, and that the tree
consisting only of the root is a multiplicative identity (see~\cite{Grossman}
or~\cite{GL:families} for details).
 
We define the comultiplication $\Delta$
on $k\{\LT(E_1$, \ldots, $E_M)\}$ as follows.
Let $t$ be a labeled tree, and let $s_1$, \ldots, $s_r$ be the children
of the root of $t$.
If $P$ is a subset of $C_t=\{s_1$, \ldots, $s_r\}$, let
$t_P$ be the labeled tree formed by making the
elements of $P$ the children of a new root, keeping the original
labels.
Define
$\Delta(t)=\sum_{P\subseteq C_t}t_P\otimes t_{C_t\backslash P}$,
where $X\backslash Y$ denotes the set-theoretic complement of
$Y$ in $X$.
Define the augmentation $\epsilon(t)$ of the bialgebra
to be $1$ if $t$ has only one node (its root),
and $0$ otherwise.
We define a grading on $k\{\LT(E_1$, \ldots, $E_M)\}$ by letting
$k\{\LT(E_1$, \ldots, $E_M)\}_n$ be the
subspace of $k\{\LT(E_1$, \ldots, $E_M)\}$ spanned by the trees
with $n+1$ nodes.
The following theorems are proved in \cite{GL:families}.
 
\begin{theorem}
$H=k\{\LT(E_1$, \ldots, $E_M)\}$ is a cocommutative graded
connected bialgebra.
\end{theorem}
 
\begin{theorem} \label{prims}
The set of labeled trees $t$ whose root has exactly one child is a
basis for the primitives $P(H)$ of $H=k\{\LT(E_1$, \ldots, $E_M)\}$.
\end{theorem}
 
Let $R$ denote the field of rational functions in the variables
$x_1,\ldots, x_N$ with coefficients from the
field $k$, and let  $E_1,\ldots, E_M$ denote $M$ derivations of
$R$ of the form
\[
E_{\gamma}=\sum_{\mu=1}^N b_{\gamma}^{\mu}
{{\partial}\over{\partial x_{\mu}}},
\]
where $b_{\gamma}^{\mu}\in R$.
We now define an $H$-module algebra on $R$.
The action of $H$ on $R$ is given by the map
$\psi: H \rightarrow \End{k}R$,
which is defined as follows.
\begin{enumerate}
\item Given a labeled tree $t$ with $m+1$ nodes,
assign the root the number
$0$ and assign the remaining nodes the numbers $1$, \ldots, $m$.
We identify the node with the number
assigned to it.  To the node $k$  asociate the summation index
$\mu_k.$  Denote $(\mu_1$, \ldots, $\mu_m)$ by $\mu$.
\item For the labeled tree $t$, let $k$ be a node of $t$,
labeled with $E_{\gamma_k}$ if $k\not=0$,
and
let $l$, \ldots, $l'$ be the children of $k$.  Define
\[
c(k;\mu) = \left\{
\begin{array}{ll}
\displaystyle
\prtl{\mu_l} \cdots \prtl{\mu_{l'}}b_{\gamma_k}^{\mu_k}(x)
& \mbox{if $k\not=0$ is not the root;} \\[7pt]
\displaystyle
\prtl{\mu_l} \cdots \prtl{\mu_{l'}}
& \mbox{if $k=0$ is the root.}
\end{array}\right.
\]
Note that if $k \neq 0$, then $c(k;\mu) \in R$.
\item Define
\[
\psi(t) = \sum_{\mu_1,\ldots, \mu_m=1}^N c(m;\mu) \cdots
c(1;\mu)c(0;\mu).
\]
\item Extend $\psi$ to all of $H$ by linearity.
\end{enumerate}
It is straight-forward to check that this action of $H$ on $R$ makes
$R$ into a $H$-module algebra.
 
An element $p\in\Hstar$ can be thought of as an infinite
series whose terms are indexed by labeled trees rather than by words,
as well as an element of a power series algebra.
Theorem~\ref{maintheorem} gives necessary and sufficient conditions
for $p$ to be differentially produced in this case.

\end{document}